\documentclass[12pt,a4paper,reqno,dvips]{article}%
\usepackage{amsfonts}
\usepackage{amsmath}
\usepackage{amssymb}
\usepackage{graphicx}
\newtheorem{theorem}{Theorem}

\newtheorem{definition}[theorem]{Definition}

\newtheorem{proposition}[theorem]{Proposition}

\title{\bf Iterated Differential Forms IV: $\mathcal{C}$--Spectral Sequence}
\author{\sc{A.~M.~Vinogradov}\thanks{{\bf e}-{\it mail}: \texttt{vinograd@unisa.it}} and \sc{L.~Vitagliano}\thanks{{\bf e}-{\it mail}: \texttt{luca\_vitagliano@fastwebnet.it}}\\
\small{DMI, Universit\`a degli Studi di Salerno}\\ \small{and INFN, Gruppo collegato di Salerno,}\\
\small{Via Ponte don Melillo, 84084 Fisciano (SA), Italy}}
\begin{document}
\maketitle
\begin{abstract}
For the multiple differential algebra of iterated differential
forms \cite{vv06,vv06b,ks03,s06} on a diffiety
$(\mathcal{O},\mathcal{C})$ \cite{b99,v01} an analogue of
$\mathcal{C}$--spectral sequence \cite{v78,v84,v01} is
constructed. The first term of it is naturally interpreted as the
algebra of \emph{secondary} \cite{b99,v01,v98} iterated
differential forms on $(\mathcal{O},\mathcal{C})$. This allows to
develop \emph{secondary tensor analysis} on generic diffieties,
some simplest elements of which are sketched here. The presented
here general theory will be specified to infinite jet spaces and
infinitely prolonged PDEs in subsequent notes.
\end{abstract}

\maketitle
\newpage
\section{Introduction}

Let $(\mathcal{O},\mathcal{C})$ be a diffiety. Recall that an
elementary diffiety is a pair
$(\mathcal{E}^{\infty},\mathcal{C}(\mathcal{E}))$, where
$\mathcal{E}^{\infty}$ is the infinite prolongation of a system
$\mathcal{E}$ of partial differential equations and
$\mathcal{C}(\mathcal{E})$ the Cartan distribution on it (see
\cite{b99,v78,v84}). More generally, a diffiety is a smooth,
usually pro--finite manifold $\mathcal{O}$ endowed with a
finite dimensional integrable distribution $\mathcal{C}$, locally
identical to an elementary diffiety (see \cite{v01,v98}).

Recall that Secondary Calculus on the diffiety
$(\mathcal{E}^{\infty},\mathcal{C}(\mathcal{E}))$ is the proper
formalization of the idea of differential calculus on the \textquotedblleft space
of all solutions of a system of partial differential equations
$\mathcal{E}$\textquotedblright. By this reason the search for analogues in
secondary calculus of objects, concepts and constructions of
\textquotedblleft traditional\textquotedblright{} mathematics is an actual and important problem,
called the \emph{secondarization problem} (see \cite{v01,v98}). In
this note this problem is solved for iterated differential forms
(in the sequel, IDFs) on the basis of the algebraic approach
developed in preceding notes \cite{vv06,vv06b} of this series.
Below we use the notation, definitions and results of these notes.
It is worth noticing that the alternative approach to iterated
forms proposed in \cite{ks03,s06} is not very secondarization
friendly.

Recall that the solution of the secondarization problem for usual,
i.e., one time iterated forms, is given in \cite{vv06,vv06b} (see
also \cite{v01,v98}). Namely, the algebra of secondary
differential forms is realized as the first term of the
$\mathcal{C}$--spectral sequence. In this note an analogue of this
sequence is constructed and its first term is interpreted as the
algebra of \emph{secondary IDFs}.

\section{Preliminaries on IDF\lowercase{s}}

Let $\mathcal{O}$ be a smooth, generally, pro--finite manifold,
i.e., the inverse limit of a sequence of finite--dimensional
submersions
\begin{equation}\label{Seq}
\cdots\longleftarrow\mathcal{O}_{s-1}\longleftarrow\mathcal{O}_{s}%
\longleftarrow\mathcal{O}_{s+1}\longleftarrow\cdots.
\end{equation}
Denote by
$(\Lambda_{k}(\mathcal{O}_{s}),d_{1},\ldots,d_{k})$ the $%
\mathbb{Z}
^{k}$--graded multiple differential algebra of geometric
$k$--times IDFs on the (usual) manifold $\mathcal{O}_{s}$ (see
\cite{vv06,vv06b}). Sequence (\ref{Seq}) induces the sequence of
embeddings of graded multiple differential algebras
\[
\cdots\subset\Lambda_{k}(\mathcal{O}_{s-1})\subset\Lambda_{k}(\mathcal{O}%
_{s})\subset\Lambda_{k}(\mathcal{O}_{s+1})\subset\cdots.
\]
Then $\Lambda_{k}(\mathcal{O})\overset{\mathrm{def}}{=}%
{\textstyle\bigcup\nolimits_{s}}
\Lambda_{k}(\mathcal{O}_{s})$ is a $%
\mathbb{Z}
^{k}$--graded multiple differential algebra.

\begin{definition}
Elements of the algebra $\Lambda_{k}(\mathcal{O})$ are called
$k$--times IDFs on the pro--finite manifold $\mathcal{O}$.
\end{definition}

The theory of IDFs on pro--finite manifolds is a straightforward
generalization of the theory presented in \cite{vv06,vv06b} and by
this reason we shall not go into its details here.

Denote by $\mathrm{D}(\Lambda_{k}(\mathcal{O}),\Lambda
_{m}(\mathcal{O}))$ the $\Lambda_{m}(\mathcal{O})$--module of
graded $\Lambda_{m}(\mathcal{O})$--valued derivations of
$\Lambda_{k}(\mathcal{O})$, $k\leq m$.

\begin{proposition}
Let $X\in\mathrm{D}(C^{\infty}(\mathcal{O}),\Lambda_{k}(\mathcal{O}))$ and
$K\subset\{1,\ldots,k\}$. Then there exists a unique $i_{X}^{K}\in\mathrm{D}%
(\Lambda_{k}(\mathcal{O}),\Lambda_{k}(\mathcal{O}))$ such that
\[
i_{X}^{K}(d_{K^{\prime}}f)=\left\{
\begin{array}
[c]{cc}%
d_{K^{\prime}\smallsetminus K}(X(f)), & \text{if }K\subset K^{\prime}\\
0, & \textrm{otherwise}%
\end{array}
\right.  ,\quad K^{\prime}\subset\{1,\ldots,k\}.
\]
\end{proposition}
For instance, if $i_{m}^{K}\overset{\mathrm{def}}{=}i_{d_{m}|_{C^{\infty}(\mathcal{O})}%
}^{K}\in\mathrm{D}(\Lambda_{k}(\mathcal{O}),\Lambda_{k}(\mathcal{O}))$,
$m\leq k$, $K\subset\{1,\ldots,k\}$, then
$i_{m}^{\varnothing}=d_{m}$.

Recall that, according to proposition 3 of \cite{vv06} covariant
$k$--tensors on $\mathcal{O}$ are naturally interpreted as
$k$--times IDFs. Namely, for
any $k$ there exists an injective morphism of $C^{\infty}(\mathcal{O}%
)$--modules
$\iota_{k}:T_{k}^{0}(\mathcal{O})\hookrightarrow\Lambda
_{k}(\mathcal{O})$. Then proposition 4 of \cite{vv06}
characterizes those IDFs that represent covariant tensors. Below
we give another useful characterization of such IDFs.

\begin{proposition}
\label{CharCovTens}A homogeneous element $\omega\in\Lambda_{k}(\mathcal{O})$
of multi--degree $(1,\ldots,1,1)\in%
\mathbb{Z}
^{k}$ is a covariant tensor (i.e.,
$\omega\in\operatorname{im}\iota_{k}$) iff
$i_{m}^{K}(\omega)=0$ for any $m< k$ and $K=\{k_{1}%
,\ldots,k_{s}\}\subset\{1,\ldots,k\}$, $s\geq2$.
\end{proposition}

\section{Secondary IDF\lowercase{s}}

Let $(\Lambda_{k}(\mathcal{O}),d_{1},\ldots,d_{k})$ be the multiple complex of
$k$--times IDFs on $\mathcal{O}$ and $\kappa_{1k}:\Lambda_{k}(\mathcal{O}%
)\longrightarrow\Lambda_{k}(\mathcal{O})$ the involution that
interchanges differentials $d_{1}$ and $d_{k}$, $k\geq1$. The
ideal $\mathcal{C}\Lambda(\mathcal{O})\subset\Lambda(\mathcal{O})$
of (ordinary) differential forms on $\mathcal{O}$ vanishing on
$\mathcal{C}(\mathcal{O})$ is embedded naturally into
$\Lambda_{k}(\mathcal{O})$ due to the embedding
$\Lambda(\mathcal{O})\equiv\Lambda_{1}(\mathcal{O})\subset\Lambda_{k}(\mathcal{O})$.
Denote by
$\mathcal{C}\Lambda_{k}(\mathcal{O})\subset\Lambda_{k}(\mathcal{O})$
the ideal generated by elements in the form
\[
(d_{K}\circ\kappa_{1k})(\omega),\quad\omega\in\mathcal{C}\Lambda
(\mathcal{O}),\quad K\subset\{1,\ldots,k-1\}.
\]
It is multi--differential, i.e.,
\begin{equation}
d_{m}(\mathcal{C}\Lambda_{k}(\mathcal{O}))\subset\mathcal{C}\Lambda
_{k}(\mathcal{O}),\quad m\leq k. \label{Comp}%
\end{equation}

\begin{definition}
Elements of $\mathcal{C}\Lambda_{k}(\mathcal{O})$ are called
$k$--Cartan $k$--times IDFs.
\end{definition}
Let $\mathcal{C}^{p}\Lambda_{k}(\mathcal{O})$ be the $p$--th power
of $\mathcal{C}\Lambda_{k}(\mathcal{O})$. Then the sequence
\begin{equation}
\Lambda_{k}(\mathcal{O})\supset\mathcal{C}\Lambda_{k}(\mathcal{O}%
)\supset\mathcal{C}^{2}\Lambda_{k}(\mathcal{O})\supset\cdots\supset
\mathcal{C}^{p}\Lambda_{k}(\mathcal{O})\supset\cdots, \label{Cfilt}%
\end{equation}
is a filtration of the multi--graded differential algebra $(\Lambda
_{k}(\mathcal{O}),d_{k})$ by means of multi--graded differential ideals.

\begin{definition}
The spectral sequence $\Lambda_{k-1}\mathcal{C}E(\mathcal{O})=
\{(\Lambda_{k-1}\mathcal{C}E_{r}(\mathcal{O}),d_{k,r})\}_{r}$
associated with
filtration (\ref{Cfilt}) is called the $\Lambda_{k-1}\mathcal{C}%
$\emph{--spectral sequence} of the diffiety $(\mathcal{O},\mathcal{C})$.
\end{definition}
In particular,
$\Lambda_{0}\mathcal{C}E(\mathcal{O})=\mathcal{C}E(\mathcal{O})$.

Filtration (\ref{Cfilt}) is regular and, therefore, the $\Lambda
_{k-1}\mathcal{C}$--spectral sequence converges to the cohomology
$H(\Lambda _{k}(\mathcal{O}),d_{k})$. Hence,
\[
\Lambda_{k-1}\mathcal{C}E_{\infty}(\mathcal{O})\simeq H(\Lambda_{k}%
(\mathcal{O}),d_{k})\simeq H(\mathcal{O}).
\]

It follows from (\ref{Comp}) and commutation relations
$[d_{m},d_{k}]=0$, $m<k$, that $d_{m}$ induces a well defined
differential in each term of the
$\Lambda_{k-1}\mathcal{C}$--spectral sequence.  Denote by
$d_{m,r}$ the so-obtained differential in the term
$\Lambda_{k-1}\mathcal{C}E_{r}(\mathcal{O})$. Obviously,
$[d_{m,r},d_{n,r}]=0$, $m,n\leq k$, so that $(\Lambda_{k-1}\mathcal{C}%
E_{r}(\mathcal{O}),d_{1,r},\ldots,d_{k,r})$ is a
multiple differential algebra.

The zeroth term of the $\Lambda_{k-1}\mathcal{C}$--spectral
sequence is described as follows. Let
$\mathcal{C}^{p}\Lambda_{k}^{p}(\mathcal{O})$ denote the
$\Lambda_{k-1}(\mathcal{O})$--module of $k$--times iterated
$p$--forms in
$\mathcal{C}^{p}\Lambda_{k}(\mathcal{O})$. Note that $\mathcal{C}^{p}%
\Lambda_{k}(\mathcal{O})=\mathcal{C}^{p}\Lambda_{k}^{p}(\mathcal{O}%
)\wedge\Lambda_{k}(\mathcal{O})$ and put
\[
\mathcal{H}\Lambda_{k}(\mathcal{O})\overset{\mathrm{def}}{=}\Lambda
_{k}(\mathcal{O})\,/\,\mathcal{C}\Lambda_{k}(\mathcal{O})=\Lambda
_{k-1}\mathcal{C}E_{0}(\mathcal{O}).
\]

\begin{definition}
$(\mathcal{H}\Lambda_{k}(\mathcal{O}),d_{k,0})$ is called the
differential algebra of \emph{horizontal }$k$\emph{--times IDFs}
on the diffiety $(\mathcal{O},\mathcal{C})$.
\end{definition}

\begin{proposition}
\label{E0}There exists a natural isomorphism
\[\Lambda_{k-1}\mathcal{C}E_{0}^{p,\bullet}(\mathcal{O}%
)\simeq\mathcal{C}^{p}\Lambda_{k}^{p}(\mathcal{O})\otimes_{\Lambda_{k-1}%
}\mathcal{H}\Lambda_{k}(\mathcal{O}).\]
\end{proposition}

The isomorphism of proposition \ref{E0} supplies $\mathcal{C}^{p}\Lambda_{k}^{p}(\mathcal{O}%
)\otimes_{\Lambda_{k-1}}\mathcal{H}\Lambda_{k}(\mathcal{O})$ with
a structure of a differential algebra. The corresponding
differential $\overline{d}_{k}$ looks as follows. Let
$\omega\in\mathcal{C}^{p}\Lambda _{k}^{p}(\mathcal{O})$ and
$\overline{\rho}=[\rho]_{\mathcal{C}\Lambda
_{k}(\mathcal{O})}\in\mathcal{H}\Lambda_{k}(\mathcal{O})$,
$\rho\in\Lambda _{k}(\mathcal{O})$. In its turn
$d_{k}\omega\in\mathcal{C}^{p}\Lambda
_{k}(\mathcal{O})$. Therefore, $d_{k}\omega=\sum_{\alpha}\omega_{\alpha}%
\wedge\sigma_{\alpha}$ for some $\omega_{\alpha}\in\mathcal{C}^{p}\Lambda
_{k}^{p}(\mathcal{O})$ and $\sigma_{\alpha}\in\Lambda_{k}(\mathcal{O})$. Then
\[
\overline{d}_{k}(\omega\otimes\overline{\rho})=%
{\textstyle\sum\nolimits_{\alpha}}
\omega_{\alpha}\otimes\overline{\sigma_{\alpha}\wedge\rho}+(-1)^p\omega\otimes
d_{k,0}\overline{\rho}\in\mathcal{C}^{p}\Lambda_{k}^{p}(\mathcal{O}%
)\otimes_{\Lambda_{k-1}}\mathcal{H}\Lambda_{k}(\mathcal{O}),
\]
where $\overline{\sigma_{\alpha}\wedge\rho}=[\sigma_{\alpha}\wedge
\rho]_{\mathcal{C}\Lambda_{k}(\mathcal{O})}\in\mathcal{H}\Lambda
_{k}(\mathcal{O})$.

\begin{definition}
\label{SecIter}$(\Lambda_{k-1}\mathcal{C}E_{1}(\mathcal{O}),d_{1,1}%
,\ldots,d_{k,1})$ is called the multiple differential algebra of
\emph{secondary} $k$\emph{--times IDFs}. Accordingly, elements of
$\Lambda_{k-1}\mathcal{C}E_{1}^{p,\bullet}(\mathcal{O})$ are
called \emph{secondary }$k$\emph{--times iterated differential
}$p$\emph{--forms}.
\end{definition}

It takes place the secondary analogue of the fact that
$(k+1)$--times iterated $0$--forms coincide with $k$--times
iterated ones. Indeed, it holds the

\begin{proposition}
\label{Embed} There exists a natural isomorphism of graded multiple
differential algebras
\[
\varphi_{k-1}:(\Lambda_{k-1}\mathcal{C}E_{1}(\mathcal{O}),d_{1,1}%
,\ldots,d_{k,1})\longrightarrow(\Lambda_{k}\mathcal{C}%
E_{1}^{0,\bullet}(\mathcal{O}),d_{1,1},\ldots,d_{k,1}).
\]

\end{proposition}

Since $\Lambda_{k-1}\mathcal{C}E_{1}^{0,\bullet}(\mathcal{O})=H(\mathcal{H}%
\Lambda_{k}(\mathcal{O}),d_{k,0})$, according to the above
proposition secondary $k$--times IDFs may be understood either as
elements in the first term of the
$\Lambda_{k-1}\mathcal{C}$--spectral sequence or as cohomologies
of the differential algebra of \emph{horizontal}
$(k+1)$\emph{--times IDFs}.

The standard action of the permutation group $S_{k}$ on IDFs is secondarized
as follows. Let $\sigma\in S_{k-1}$ be a permutation of $\{1,\ldots,k-1\}$.
The associated automorphism $\kappa_{\sigma}:\Lambda_{k}(\mathcal{O}%
)\longrightarrow\Lambda_{k}(\mathcal{O})$ respects both the
differential $d_{k}$ and the filtration (\ref{Cfilt}), i.e.,
$\kappa_{\sigma}\circ d_{k}=d_{k}\circ\kappa_{\sigma}$ and
$\kappa_{\sigma}(\mathcal{C}\Lambda
_{k}(\mathcal{O}))\subset\mathcal{C}\Lambda_{k}(\mathcal{O})$.
Therefore,
$\kappa_{\sigma}$ induces an automorphism of the $\Lambda_{k-1}\mathcal{C}%
$--spectral sequence. Abusing the notation we continue denoting by
$\kappa _{\sigma}$ the induced automorphism in the first term of
$\Lambda
_{k-1}\mathcal{C}E(\mathcal{O})$. In particular, $\kappa_{\sigma}:\Lambda_{k-1}%
\mathcal{C}E_{1}(\mathcal{O})\longrightarrow\Lambda_{k-1}\mathcal{C}%
E_{1}(\mathcal{O})$ is an isomorphism of the multiple differential
algebras
$(\Lambda_{k-1}\mathcal{C}E_{1}(\mathcal{O}),d_{1,1},\ldots,d_{k,1})$
and
$(\Lambda_{k-1}\mathcal{C}E_{1}(\mathcal{O}),d_{\sigma(1),1},\ldots
,d_{\sigma(k-1),1},d_{k,1})$. Proposition \ref{Embed} allows to
define the action of the transposition $\tau$ exchanging $k$ and
$m<k$ on $\Lambda_{k-1}\mathcal{C}%
E_{1}(\mathcal{O})$. Namely, observe that the action of
$\kappa_{\tau}$ on
$\Lambda_{k}\mathcal{C}E_{1}^{0,\bullet}(\mathcal{O})$ is well
defined as it was explained above and the isomorphism
$\varphi_{k-1}$ transfers it then to
$\Lambda_{k-1}\mathcal{C}E_{1}(\mathcal{O})$.

Since $\mathcal{C}\Lambda_{k}(\mathcal{O})$ is a multi--graded
differential ideal,
$\Lambda_{k-1}\mathcal{C}E_{1}^{p,\bullet}(\mathcal{O})$ inherits
a
multi--graded algebra structure. Namely, let $\theta=[[\omega]_{\mathcal{C}%
^{p+1}\Lambda_{k}(\mathcal{O})}]_{\mathrm{\operatorname{im}}d_{k,0}}\in
\Lambda_{k-1}\mathcal{C}E_{1}^{p,\bullet}(\mathcal{O})$, $\omega\in
\mathcal{C}^{p}\Lambda_{k}(\mathcal{O})$, $d_{k}\omega\in\mathcal{C}%
^{p+1}\Lambda_{k}(\mathcal{O})$. Then
\begin{equation}
\theta\in\Lambda_{k-1}\mathcal{C}E_{1}^{(p_{1},\ldots,p_{k-1},p),q}%
(\mathcal{O})\overset{\mathrm{def}}{\Longleftrightarrow}\omega\text{ has
multi--degree }(p_{1},\ldots,p_{k-1},p+q)\in%
\mathbb{Z}
^{k}, \label{multigr}%
\end{equation}
$p_{1},\ldots,p_{k-1},q\geq0$.  Definition (\ref{multigr}) is
correct and
\[
\Lambda_{k-1}\mathcal{C}E_{1}^{p,\bullet}(\mathcal{O})=\bigoplus_{p_{1}%
,\ldots,p_{k-1},q}\Lambda_{k-1}\mathcal{C}E_{1}^{(p_{1},\ldots,p_{k-1}%
,p),q}(\mathcal{O}).
\]

\begin{proposition}
Let $\sigma\in S_{k}$ be a permutation of \;$\{1,\ldots,k\}$. Then
\[
\kappa_{\sigma^{-1}}(\Lambda_{k-1}\mathcal{C}E_{1}^{(p_{1},\ldots,p_{k-1}%
,p_{k}),q})\subset\Lambda_{k-1}\mathcal{C}E_{1}^{(p_{\sigma(1)},\ldots
,p_{\sigma(k-1)},p_{\sigma(k)}),q}.
\]
\end{proposition}

\section{IDF--symmetries of a diffiety}

Iterated analogues of contact and, respectively, trivial contact
vector fields on a diffiety are defined as follows (see
\cite{b99,v01}):
\[
\mathrm{D}_{\mathcal{C}}(\Lambda_{k-1}(\mathcal{O}))\overset{\mathrm{def}}%
{=}\{X\in\mathrm{D}(\Lambda_{k-1}(\mathcal{O}),\Lambda_{k-1}(\mathcal{O}%
))\;|\;\mathcal{L}_{X}^{\{k\}}(\mathcal{C}\Lambda_{k}(\mathcal{O}%
))\subset\mathcal{C}\Lambda_{k}(\mathcal{O})\}
\]
and
\[
\mathcal{C}\mathrm{D}(\Lambda_{k-1}(\mathcal{O}))\overset{\mathrm{def}}%
{=}\{X\in\mathrm{D}(\Lambda_{k-1}(\mathcal{O}),\Lambda_{k-1}(\mathcal{O}%
))\;|\;i_{X}^{\{k\}}(\mathcal{C}\Lambda_{k}(\mathcal{O}))=0\}.
\]
$\mathrm{D}_{\mathcal{C}}(\Lambda_{k-1}(\mathcal{O}))$ is a graded
sub--algebra of the Lie algebra
$\mathrm{D}(\Lambda_{k-1}(\mathcal{O}))$ and
$\mathcal{C}\mathrm{D}(\Lambda_{k-1}(\mathcal{O}))$ is an its
graded ideal. Denote by
\[
\Lambda_{k-1}\mathrm{Sym}(\mathcal{O})\overset{\mathrm{def}}{=}\mathrm{D}%
_{\mathcal{C}}(\Lambda_{k-1}(\mathcal{O}))\,/\,\mathcal{C}\mathrm{D}%
(\Lambda_{k-1}(\mathcal{O})).
\]
the quotient Lie algebra. In particular, $\Lambda_{0}\mathrm{Sym}%
(\mathcal{O})=\mathrm{Sym}(\mathcal{O})$.

\begin{definition}
Elements in $\Lambda_{k-1}\mathrm{Sym}(\mathcal{O})$ are called
$(k-1)$--IDF--symmetries of the diffiety
$(\mathcal{O},\mathcal{C})$.
\end{definition}

Just as symmetries of $(\mathcal{O},\mathcal{C})$ act naturally on
secondary differential forms, IDF--symmetries of
$(\mathcal{O},\mathcal{C})$ act
naturally on secondary IDFs. Namely, let $\chi=[X]\in\Lambda_{k-1}%
\mathrm{Sym}(\mathcal{O})$, $X\in\mathrm{D}_{\mathcal{C}}(\Lambda
_{k-1}(\mathcal{O}))$. The Lie derivative $\mathcal{L}_{X}^{\{k\}}:\Lambda
_{k}(\mathcal{O})\longrightarrow\Lambda_{k}(\mathcal{O})$ respects both the
differential $d_{k}$ and the filtration (\ref{Cfilt}), i.e., $\mathcal{L}%
_{X}^{\{k\}}\circ d_{k}=d_{k}\circ\mathcal{L}_{X}^{\{k\}}$ and $\mathcal{L}%
_{X}^{\{k\}}(\mathcal{C}\Lambda_{k}(\mathcal{O}))\subset\mathcal{C}\Lambda
_{k}(\mathcal{O})$. Therefore, $\mathcal{L}_{X}^{\{k\}}$ induces a
morphism of the $\Lambda_{k-1}\mathcal{C}$--spectral sequence. Denote by
\[
\mathcal{L}_{\chi}^{\{k\}}:\Lambda_{k-1}\mathcal{C}E_{1}^{p,\bullet
}(\mathcal{O})\longrightarrow\Lambda_{k-1}\mathcal{C}E_{1}^{p,\bullet
}(\mathcal{O})
\]
its action on the first term since it doesn't depend on the choice
of the representative $X$ of $\chi$. $\mathcal{L}_{\chi}^{\{k\}}$
is a graded derivation of
$\Lambda_{k-1}\mathcal{C}E_{1}(\mathcal{O})$. Moreover, if
$\chi_{1},\chi_{2}\in\Lambda_{k-1}\mathrm{Sym}(\mathcal{O})$, then
\[
\lbrack\mathcal{L}_{\chi_{1}}^{\{k\}},d_{k,1}]=0,\quad\lbrack\mathcal{L}%
_{\chi_{1}}^{\{k\}},\mathcal{L}_{\chi_{2}}^{\{k\}}]=\mathcal{L}_{[\chi
_{1},\chi_{2}]}^{\{k\}}.
\]

The definition of $\mathcal{L}_{\chi}^{\{k\}}$ is consistent with
the natural isomorphism
$\Lambda_{k-1}\mathcal{C}E_{1}(\mathcal{O})\simeq\Lambda
_{k}\mathcal{C}E_{1}^{0,\bullet}(\mathcal{O})$. Namely, let $\chi$
and $X$ be as above, so that
$\mathcal{L}_{X}^{\{k\}}\in\mathrm{D}_{\mathcal{C}}(\Lambda
_{k}(\mathcal{O}))$. Put $\chi^{\prime}\overset{\mathrm{def}}{=}%
[\mathcal{L}_{X}^{\{k\}}]\in\Lambda_{k}\mathrm{Sym}(\mathcal{O})$. Then
\[
\mathcal{L}_{\chi}^{\{k\}}=\varphi_{k-1}^{-1}\circ\mathcal{L}_{\chi^{\prime}%
}^{\{k+1\}}|_{\Lambda_{k}\mathcal{C}E_{1}^{0,\bullet}(\mathcal{O})}%
\circ\varphi_{k-1}.
\]

Now observe that
$i_{m}^{K}\in\mathrm{D}_{\mathcal{C}}(\Lambda_{k-1}(\mathcal{O}))$
for any $m<k$ and $K\subset\{1,\ldots,k-1\}$ and put
\[
I_{m}^{K}\overset{\mathrm{def}}{=}[i_{m}^{K}]\in\Lambda_{k-1}\mathrm{Sym}%
(\mathcal{O}).
\]
Then $d_{m,1}=\mathcal{L}_{I_{m}^{\varnothing}}^{\{k\}}$.

IDF--symmetries of $(\mathcal{O},\mathcal{C})$ can be also
inserted into secondary IDFs as follows. Indeed, let $\chi$ and $X$ be as
above and $\theta=[[\omega
]_{\mathcal{C}^{p+1}\Lambda_{k}(\mathcal{O})}]_{\mathrm{\operatorname{im}%
}\;d_{k,0}}\in\Lambda_{k-1}\mathcal{C}E_{1}^{p,\bullet}(\mathcal{O})$,
$\omega\in\mathcal{C}^{p}\Lambda_{k}(\mathcal{O})$,
$d_{k}\omega\in \mathcal{C}^{p+1}\Lambda_{k}(\mathcal{O})$, be a
secondary IDF. Then
$i_{X}^{\{k\}}\omega\in\mathcal{C}^{p-1}\Lambda_{k}(\mathcal{O})$
and
$d_{k}(i_{X}^{\{k\}}\omega)=\mathcal{L}_{X}^{\{k\}}\omega-i_{X}^{\{k\}}%
(d_{k}\omega)\in\mathcal{C}^{p}\Lambda_{k}(\mathcal{O})$ so that
\[
i_{\chi}^{\{k\}}\theta\overset{\mathrm{def}}{=}[[i_{X}^{\{k\}}\omega
]_{\mathcal{C}^{p}\Lambda_{k}(\mathcal{O})}]_{\mathrm{\operatorname{im}%
}\;d_{k,0}}\in\Lambda_{k-1}\mathcal{C}E_{1}^{p-1,\bullet}(\mathcal{O})
\]
is a well defined secondary IDF. Moreover, the map
\[
i_{\chi}^{\{k\}}:\Lambda_{k-1}\mathcal{C}E_{1}^{p,\bullet}(\mathcal{O}%
)\longrightarrow\Lambda_{k-1}\mathcal{C}E_{1}^{p-1,\bullet}(\mathcal{O})
\]
is a graded derivation of
$\Lambda_{k-1}\mathcal{C}E_{1}(\mathcal{O})$ not
depending on the choice of $X$ in $\chi$.

For $\chi_{1},\chi_{2}%
\in\Lambda_{k-1}\mathrm{Sym}(\mathcal{O})$ the following secondary
analogues of well-known (graded) commutation relations hold
\[
\lbrack i_{\chi_{1}}^{\{k\}},d_{k,1}]=\mathcal{L}_{\chi_{1}}^{\{k\}}%
,\quad\lbrack i_{\chi_{1}}^{\{k\}},i_{\chi_{2}}^{\{k\}}]=0,\quad\lbrack
i_{\chi_{1}}^{\{k\}},\mathcal{L}_{\chi_{2}}^{\{k\}}]=i_{[\chi_{1},\chi_{2}%
]}^{\{k\}}.
\]
The definition of $i_{\chi}^{\{k\}}$ is consistent with the
natural isomorphism
$\Lambda_{k-1}\mathcal{C}E_{1}(\mathcal{O})\simeq\Lambda
_{k}\mathcal{C}E_{1}^{0,\bullet}(\mathcal{O})$. Namely, if $\chi$
and $X$ be
as above, then $i_{X}^{\{k\}}\in\mathrm{D}_{\mathcal{C}}(\Lambda_{k}(\mathcal{O}%
))$ and
\[
i_{\chi}^{\{k\}}=\varphi_{k-1}^{-1}\circ\mathcal{L}_{\chi^{\prime\prime}%
}^{\{k+1\}}|_{\Lambda_{k}\mathcal{C}E_{1}^{0,\bullet}(\mathcal{O})}%
\circ\varphi_{k-1}.
\]
with $\chi^{\prime\prime}\overset{\mathrm{def}}{=}[i_{X}^{\{k\}}%
]\in\Lambda_{k}\mathrm{Sym}(\mathcal{O})$.

\begin{definition}
\label{DefSecCovTens}Let
$\tau\in\Lambda_{k-1}\mathcal{C}E_{1}^{(1,\ldots
,1,1),\bullet}(\mathcal{O})$ be a secondary IDF of multi--degree
$(1,\ldots,1,1)$. $\tau$ is called a \emph{secondary covariant
}$k\emph{--tensor}$ iff
$\mathcal{L}_{I_{m}^{K}}^{\{k\}}\tau=i_{I_{m}^{K^{\prime}}}^{\{k\}}\tau=0$
for
any $m<k$ and $K=\{k_{1},\ldots,k_{s}\},K^{\prime}=\{k_{1}^{\prime}%
,\ldots,k_{s^{\prime}}^{\prime}\}\subset\{1,\ldots,k-1\}$, $s\geq2$,
$s^{\prime}\geq1$.
\end{definition}

Definition \ref{DefSecCovTens} mimics characterization of
covariant tensors given in proposition \ref{CharCovTens}.

\section{Other related $\mathcal{C}$--spectral sequences}

Now on put $\mathcal{C}_{k}\Lambda_{k}(\mathcal{O})=\mathcal{C}%
\Lambda_{k}(\mathcal{O})$ and observe that
$\mathcal{C}_{k}\Lambda_{k}(\mathcal{O})$ is not unique canonical
multi--differential ideal in $\Lambda_{k}(\mathcal{O})$. Indeed,
such are ideals
\[
\mathcal{C}_{m}\Lambda_{k}(\mathcal{O})\overset{\mathrm{def}}{=}\kappa
_{mk}(\mathcal{C}_{k}\Lambda_{k}(\mathcal{O})),
\]
for any $m\leq k$. Here
$\kappa_{mk}:\Lambda_{k}(\mathcal{O})\longrightarrow\Lambda_{k}(\mathcal{O})$
is the involution that interchanges differentials $d_{m}$ and
$d_{k}$.
More generally, let $K=\{m_{1}%
,\ldots,m_{r}\}\subset\{1,\ldots,k\}$. Put
\[
\mathcal{C}_{K}\Lambda_{k}(\mathcal{O})\overset{\mathrm{def}}{=}%
\mathcal{C}_{m_{1}}\Lambda_{k}(\mathcal{O})+\cdots+\mathcal{C}_{m_{r}}%
\Lambda_{k}(\mathcal{O}).
\]
Then
$\mathcal{C}_{K}\Lambda_{k}(\mathcal{O})\subset\Lambda_{k}(\mathcal{O})$
is a multi--differential for any $K\subset\{1,\ldots,k\}.$

\begin{definition}
Elements in $\mathcal{C}_{K}\Lambda_{k}(\mathcal{O})$ are called $K$--Cartan
$k$--times IDFs.
\end{definition}

All these ideals transform $(\Lambda_{k}(\mathcal{O}),d_{k})$ into
a \emph{multi--filtered} complex. Homological algebra of such kind
complexes and appearing this way numerous
$\mathcal{C}$--spectral--like sequences will be considered
separately. They enrich noteworthy secondary calculus and, as a
consequence, geometrical theory of nonlinear PDEs with new
powerful instruments.

\end{document}